# Optimal Distributed Voltage Control via Primal Dual Gradient Dynamics

Mohammed N. Khamees, *Student Member*, IEEE, Yang Liu, *Member*, *IEEE*, Kai Sun, *Senior Member*, IEEE
Department of Electrical Engineering and Computer Science, University of Tennessee, Knoxville, TN, USA
mkhamees@vols.utk.edu, yliu161@vols.utk.edu, kaisun@utk.edu

*Abstract*—The rapidly increasing penetration of inverter-based resources into a power transmission network requires more sophisticated voltage control strategies considering their inherent output variabilities. In addition, faults and load variations affect the voltage profile over the power network. This paper proposes a Primal Dual Gradient Dynamics based optimal distributed voltage control approach that optimizes outputs of distributed reactive power sources to maintain an acceptable voltage profile while preserving operational limits. Case studies of this new approach on IEEE test systems have verified its effectiveness.

*Keywords— distributed voltage control, reactive power compensation, optimization.*

## I. INTRODUCTION

The operating condition of a power transmission network continuously changes as a result of load variations, disturbances and other operational uncertainties, leading to voltage fluctuation across all buses of the network. So, facing such traditional challenges and the emerging challenges due to large-scale wind and solar integration, existing voltage control strategies need to be improved. There have been works addressing new challenges in power distribution systems but changes on existing voltage control strategies for transmission systems are more difficult [1].

Power systems have extensively implemented centralized schemes for many optimization and control functions, in which a central controller collects measurements, carries out computations and issues new commands of control [2]. In such a scheme, all the buses are required to communicate with a central controller when communication difficulties like delays, limited bandwidth, node failures, etc. are present [3]. Comparatively, a distributed control scheme employs controllers at multiple buses spreading out through the entire network, in which communication is limited to a set of neighbors of each bus. Therefore, less communication infrastructure, enhanced cybersecurity, robustness to control failure, and the ability to apply parallel computing are various potential benefits of distributed algorithms versus centralized ones [2].

A common approach to implement voltage control is applying optimal power flow with constraints imposed on bus voltages and reactive powers. After solving the problem, the reactive power injections at buses are adjusted. Therefore, this approach is using a feedforward optimization scheme, where the disturbance is presumed to be explicitly known for the controller. In contrast, there is no need to assume this explicit knowledge in feedback optimization, where the controller act based on measured data.

Although a power system has a wide range of acceptable operating conditions in which voltages and reactive powers can vary while still satisfying the limits, some of these conditions are better than the others if we consider the operational cost, transmission loss and the loss of opportunity. An example of the last is related to inverter-based renewable energy sources. When their inverters are expected to provide more reactive powers, a price is yielded since it is more economical for them to inject more real powers than reactive powers.

In this work, we exploit the feedback optimization and the distributed control strategy to propose continuous optimal feedback voltage controller, in which each bus locally communicate with its neighbors that are physically connected to share local measurements enabling the controller to adjust its reactive power output using the Primal Dual Gradient Dynamics (PDGD) algorithm. Ref. [3] has applied PDGD to optimal distributed feedback voltage control for power distribution systems, in which the radial structure of a distribution network is utilized to simplify the model of power flows for control. This paper will focus on a power transmission network having a meshed structure. The PDGD is applied in the proposed optimal distributed voltage controller to exploit the network sparsity to optimize the problem completely in distributed fashion. The controller will (1) minimize the operational cost, (2) keep voltages in acceptable ranges, and (3) satisfy reactive power limits.

The use of PDGD is motivated by the structure of the optimization problems, which provides the Lagrangian saddle-point dynamics an approach to be solved in distributed fashion. Hence, the dynamics is common in network optimization [4]. In [5], the authors use the PDGD distributed optimizer to solve primary frequency regulation load control problem. In [6], an energy-based approach is presented to study the stability of power systems coupled with market dynamics, the PDGD is used to form distributed dynamic optimization algorithm.

The rest of the paper is organized as follows: Section II introduces the power system model and provides a detailed discussion of the proposed optimal distributed feedback voltage controller. The case studies on IEEE benchmark systems are presented in section III. Finally, conclusions are drawn in section IV.

## II. PROBLEM FORMULATION

### A. Bus Injection Mode

In this paper, we consider bus injection model to represent the transmission network. Where power flow equations of an N

bus system can be written as follows in the polar form:

$$P_k = V_k \sum_{n=1}^{N} Y_{kn} V_n \cos(\delta_k - \delta_n - \theta_{kn}) \quad (1)$$

$$Q_k = V_k \sum_{n=1}^{N} Y_{kn} V_n \sin(\delta_k - \delta_n - \theta_{kn}) \quad (2)$$

where $P_k$ and $Q_k$ are the real and reactive powers injected into bus $k$. $V_k$ and $\delta_k$ are the voltage magnitude and phase angle at bus $k$. $Y_{kn}$ and $\theta_{kn}$ are the magnitude and phase angle of each element of admittance matrix. To linearize these equations, the Taylor's series expansion of a multivariable function can be applied to result in first order linear equations:

$$\begin{pmatrix} \frac{\partial \mathbf{P}}{\partial \boldsymbol{\delta}} & \frac{\partial \mathbf{P}}{\partial \mathbf{V}} \\ \frac{\partial \mathbf{P}}{\partial \boldsymbol{\delta}} & \frac{\partial \mathbf{Q}}{\partial \mathbf{V}} \end{pmatrix} \begin{pmatrix} \Delta \boldsymbol{\delta} \\ \Delta \mathbf{V} \end{pmatrix} = \begin{pmatrix} \Delta \mathbf{P} \\ \Delta \mathbf{Q} \end{pmatrix} \quad (3)$$

where the power mismatch on the right-hand side is approximated by the Jacobian matrix multiplied by the deviation of the state vector. It is clear the Jacobian matrix is not constant and to avoid some computations, one can use some assumptions which are supported by the physics of power flows in transmission lines to simplify the Jacobian matrix. Therefore, in a power transmission system that is properly designed and operated the following holds [7]:

1. Angular differences among buses are very small.
2. The line susceptance's are much larger than the line conductance's.
3. The power injected into bus is much less than which would flow if all lines from that bus were short circuited to the reference.

Therefore, applying these approximations to simplify the elements of the Jacobian matrix as follows:

$$\begin{pmatrix} -\mathbf{B} & \mathbf{G} \\ -\mathbf{G} & -\mathbf{B} \end{pmatrix} \begin{pmatrix} \Delta \boldsymbol{\delta} \\ \Delta \mathbf{V} \end{pmatrix} = \begin{pmatrix} \Delta \mathbf{P} \\ \Delta \mathbf{Q} \end{pmatrix} \quad (4)$$

where $\mathbf{G}$ and $\mathbf{B}$ are the conductance and susceptance matrices, respectively. In [3], a relaxed branch model of mesh network, which lacks the consistency in the voltage angles, is used to represent distribution network.

### B. Optimal Distributed Feedback Voltage Controller

As mentioned earlier, we build the controller to have a feedback control loop that uses the local voltage measurement at each bus with other information as input to the controller to determine the output, the injected reactive powers at time $t$. Let $\mathbf{Q}(t)$ be a given reactive power injections at instant $t$. These injections will determine the voltage profile $\mathbf{v}(t)$. Next, using the voltage profile with other information, the controller computes $\mathbf{Q}(t+1)$, the reactive powers at time $t+1$. In this input-output relationship, the controller does not need to know details behind the system, thanks to the feedback control scheme. Therefore, the controller will be injecting reactive powers at each time $t$ to control the voltages while has no control over real power, i.e., it is constant and $\Delta \mathbf{P}$ is zero. Hence, we can use (4) to find how voltages depends on reactive powers as follows:

$$\Delta \mathbf{V}(\Delta \mathbf{Q}) = \left[ -(\mathbf{G}\mathbf{B}^{-1}\mathbf{G} + \mathbf{B}) \right]^{-1} \Delta \mathbf{Q} \quad (5)$$

The controller aims to keep the voltage within the acceptable range, while satisfying the reactive power constraints, and drive the system to the optimal operating point that has the least operational cost using local voltage measurement and shared variables among neighbors. To give more room for DERs to generate more real power we consider the objective function to be the injected reactive powers that need to be minimized. Hence, the problem can be formulated as the following:

$$\min_{\mathbf{Q}} \quad f(\mathbf{Q}) = \sum_{i=1}^{C} Q_i^2 \quad (6)$$

$$\underline{v_i} < v_i(\mathbf{Q}) < \overline{v_i} \quad (7)$$

$$\underline{Q_i} < Q_i < \overline{Q_i} \quad (8)$$

where $\mathbf{Q}$ is the reactive power injections vector and it is the decision variable. $C$ is the number of controllers. $\underline{v_i}$ and $\overline{v_i}$ are the lower and upper voltage limits, respectively. $\underline{Q_i}$ and $\overline{Q_i}$ are the lower and upper reactive power injected by controller $i$. Therefore, the Lagrangian function for this optimization problem can be written as:

$$L(\mathbf{Q}, \boldsymbol{\lambda}, \boldsymbol{\mu}) = f(\mathbf{Q}) + \underline{\boldsymbol{\lambda}}^T \left( \underline{\mathbf{v}} - \mathbf{v}(\mathbf{Q}) \right) + \overline{\boldsymbol{\lambda}}^T (\mathbf{v}(\mathbf{Q}) - \overline{\mathbf{v}}) + \underline{\boldsymbol{\mu}}^T \left( \underline{\mathbf{Q}} - \mathbf{Q} \right) + \overline{\boldsymbol{\mu}}^T (\mathbf{Q} - \overline{\mathbf{Q}}) \quad (9)$$

where $\underline{\boldsymbol{\lambda}}$ and $\overline{\boldsymbol{\lambda}}$ are the Lagrangian multipliers vectors for voltage lower and upper limits, respectively. Each has a dimension equal or less than the number of load buses $M$, i.e., buses provided with control component. $\underline{\boldsymbol{\mu}}$ and $\overline{\boldsymbol{\mu}}$ are the Lagrangian multipliers vectors for reactive power injections lower and upper limits both with dimensions of $C$. All these Lagrangian multipliers act as constraints violation level indicators for the constrained variables. It is worth to note that not all load buses may have control component. Ref. [3] has used the augmented Lagrangian in which no explicit constraints is applied on the reactive power injections. Therefore, a soft thresholding function is employed with projection of reactive power injections onto constraints which caused inconsistency in updating the Lagrangian multiplier corresponding to the reactive power injections.

To solve the optimization problem the PDGD, more details can be found in [8], is employed. At each time step $t$, the measured voltage at each node is employed to update the optimization's variables. The controller performs a gradient descent for the gradient of Lagrangian with respect to $\mathbf{Q}$, the primal variable. At the same time, it calculates gradient ascent along the Lagrangian with respect to dual variables, i.e., $\boldsymbol{\lambda}$ and $\boldsymbol{\mu}$. Then, the variables are updated. Finally, the controller injects the updated values of the reactive power into the grid. Therefore, it is looking for the saddle point of this dynamical system (10).

$$\frac{\partial Q_i(t)}{\partial t} = -\left[\frac{\partial f(\mathbf{Q}(t))}{\partial Q_i(t)} + \sum_{j=1}^{M}\frac{\partial v_j(\mathbf{Q}(t))}{\partial Q_i(t)}\left(\bar{\lambda}_j(t) - \underline{\lambda}_j(t)\right) + \bar{\mu}_j(t) - \underline{\mu}_j(t)\right] \quad (10.a)$$

$$\frac{\partial \bar{\lambda}_i(t)}{\partial t} = \left[v_i(\mathbf{Q}) - \bar{v}_i\right]^+_{\bar{\lambda}_i(t)} \quad (10.b)$$

$$\frac{\partial \underline{\lambda}_i(t)}{\partial t} = \left[\underline{v}_i - v_i(\mathbf{Q})\right]^+_{\underline{\lambda}_i(t)} \quad (10.c)$$

$$\frac{\partial \bar{\mu}_i(t)}{\partial t} = \left[Q_i(t) - \bar{Q}_i\right]^+_{\bar{\mu}_i(t)} \quad (10.d)$$

$$\frac{\partial \underline{\mu}_i(t)}{\partial t} = \left[\underline{Q}_i - Q_i(t)\right]^+_{\underline{\mu}_i(t)} \quad (10.e)$$

The first term in (10.a) is the objective function partial derivative with respect to injected reactive power, equals $2Q_i(t)$ in our case. In addition, the second term contains partial derivative of voltage with respect to injected reactive power and can be calculated using (5), which is the coefficient of $\Delta\mathbf{Q}$. The voltage appearing in equations (10.b-c) is equivalent to the measured voltage. Now, we have hybrid automaton system corresponding to our dynamical system (10), due to the use of positive projection in (10.d-e). The positive projection is employed to keep the Lagrangian multipliers evaluation positive.

### III. CASE STUDY

The proposed approach is tested on two IEEE systems available in the test case archive [9] where the power flow equations are solved using MATPOWER [10]. We assume that reactive power sources are available at all load buses and can supply or consume a specified amount of reactive power. Two general cases are examined, static load and varying load with 100 MVA as the base for all cases. We use Matlab's ODE solver ode23t for all simulations.

#### A. Static Load

First, we consider load not varying with time. Starting with the IEEE 14-bus system. It has bus 1 as the slack bus. Buses 2, 3, 6 and 8 are voltage-controlled buses (PV buses) having magnitudes voltage fixed, while all others are load buses. The system is heavily loaded resulting in low voltage profile as shown in TABLE I. Each load bus provided with control component which has a specific capacity to consume or supply Q according to its physical rating. Control components in this study can supply or consume up to 20 MVar. The voltage tolerance is set to 5% for all the cases.

We run the controller in this condition and the simulation results are shown in Fig.1. It shows the reactive power injections to the left and voltage profile to the right. It illustrates that the controller is simultaneously able to bring the voltages to the acceptable range and satisfy the reactive power constraints at the minimum cost. TABLE II shows the voltage profile after convergence and Fig. 2. shows the evaluation of the cost function as the summation of all reactive power injected squared at each time step.

TABLE I.  IEEE 14-BUS SYSYTEM VOLTAGE PROFILE (NO CONTROLLER)

| $|v_1|$ | $|v_2|$ | $|v_3|$ | $|v_4|$ | $|v_5|$ | $|v_6|$ | $|v_7|$ |
|---|---|---|---|---|---|---|
| 1.0600 | 1.0450 | 1.0100 | 0.9382 | 0.9393 | 1.0700 | 0.9806 |
| $|v_8|$ | $|v_9|$ | $|v_{10}|$ | $|v_{11}|$ | $|v_{12}|$ | $|v_{13}|$ | $|v_{14}|$ |
| 1.0900 | 0.9362 | 0.9348 | .9899 | 1.0167 | 0.9927 | 0.8970 |

TABLE II. IEEE 14-BUS SYSYTEM VOLTAGE PROFILE (WITH CONTROLLER)

| $|v_1|$ | $|v_2|$ | $|v_3|$ | $|v_4|$ | $|v_5|$ | $|v_6|$ | $|v_7|$ |
|---|---|---|---|---|---|---|
| 1.0600 | 1.0450 | 1.0100 | 0.9500 | 0.9500 | 1.0700 | 1.0131 |
| $|v_8|$ | $|v_9|$ | $|v_{10}|$ | $|v_{11}|$ | $|v_{12}|$ | $|v_{13}|$ | $|v_{14}|$ |
| 1.0900 | 0.9848 | 0.9851 | 1.0209 | 1.0500 | 1.0125 | 0.9500 |

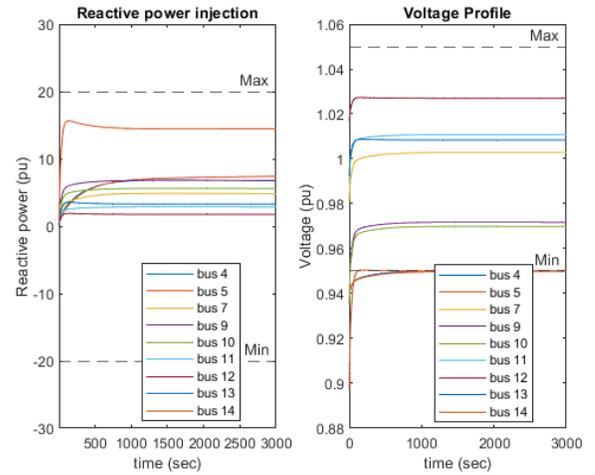

Figure 1.  IEEE 14-bus system simulation results.

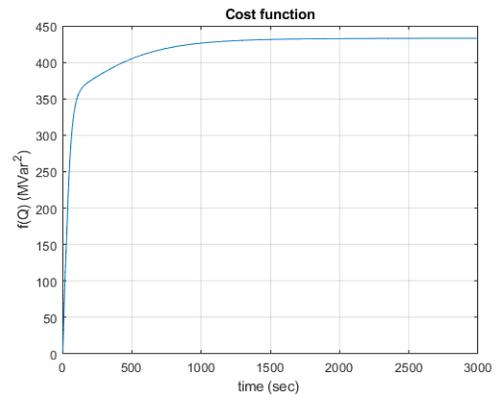

Figure 2.  IEEE 14-bus system: Cost function.

Next, we run the controller starting from the loading condition in the previous case but assume a ground fault on the transmission line between buses 4 and 5, the fault is cleared by tripping out the line. The simulation result is shown in Fig 3. Compared to the no-fault case more reactive power injected at bus 4, and the cost function has increased by 27%. In addition, we study IEEE 30-bus system which has slack bus at bus 1 and voltage-controlled buses at 2, 5, 8, 11 and 13. All others are load buses with installed control components that can supply or consume 20 MVar. This time the controller is tested under light loading to have high voltage profile. The reactive power injections and voltages evaluation are shown in Fig. 4.

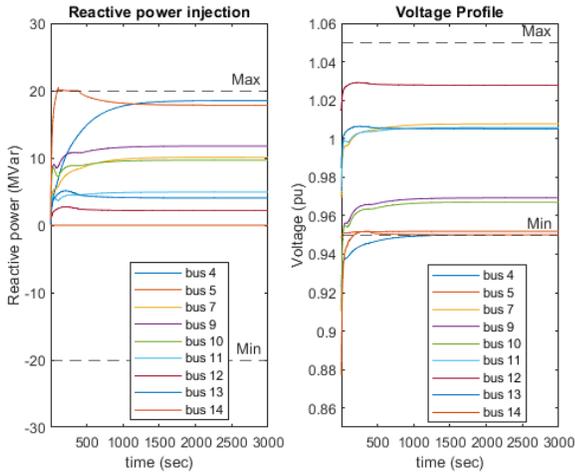

Figure 3.  IEEE 14-bus system: Line between buses 4-5 is tripped out.

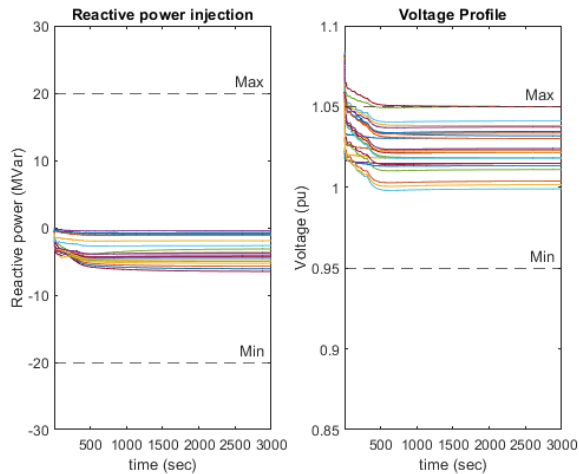

Figure 4.  IEEE 30-bus system simulation results.

## B. Varying Load

To test the controller under more practical operating conditions, we used a varying load with 24-hour time span and 1-hour time resolution, shown in Fig. 5. For the IEEE 14-bus system this load profile exists at all buses except the slack, 7, and 8 which originally have no load. So, the system is heavily loaded, causing low voltage profile at the load peak. The controller is adjusted to run every 1 hour since the operating condition has 1-hour resolution. To assess the performance of the controller the voltage profile without controller is shown in Fig. 6 and in Fig. 7 after applying the controller. The simulation result shows the controller quickly drives the voltage profile back to the acceptable range whenever it is out while keeps reactive power injection in its limits as shown in Fig 8. For the IEEE 30-bus system the voltage profile without controller highly fluctuates throughout the day as shown in Fig. 9. We run the controller under these conditions, the voltage profile for load buses is given in Fig. 10. It demonstrates that regardless of the load variations, the controller drives the voltages back into the acceptable range with the least cost. Fig. 11 shows the reactive powers injected by control components that exist over the system. It illustrates the control components consuming reactive power up to the $8^{th}$ hour then supplying reactive power to keep the voltage profile within the acceptable range.

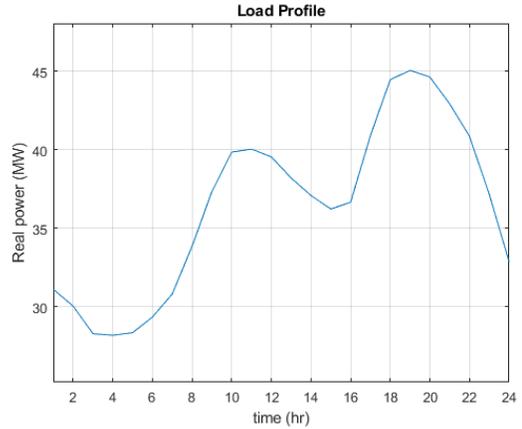

Figure 5.  IEEE 14-bus system: Load profile at all load buses.

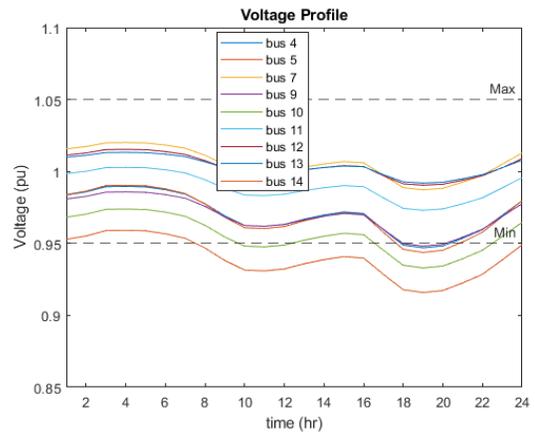

Figure 6.  IEEE 14-bus system: Voltage profile (no controller)

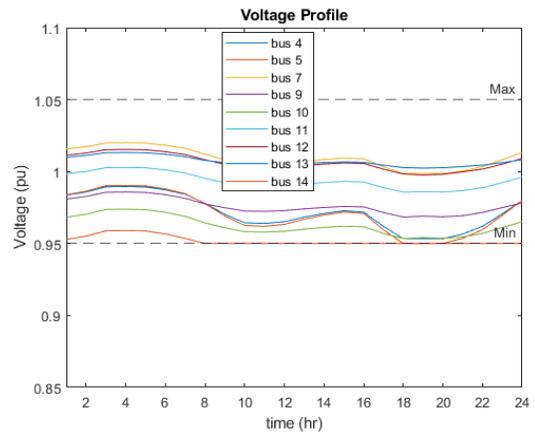

Figure 7.  IEEE 14-bus system: Voltage profile (with controller)

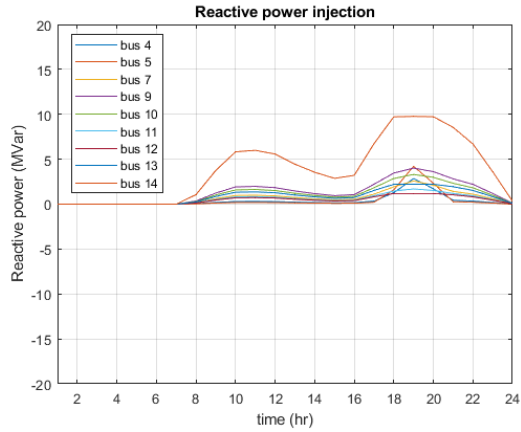

Figure 8. IEEE 14-bus system: Reactive power injection.

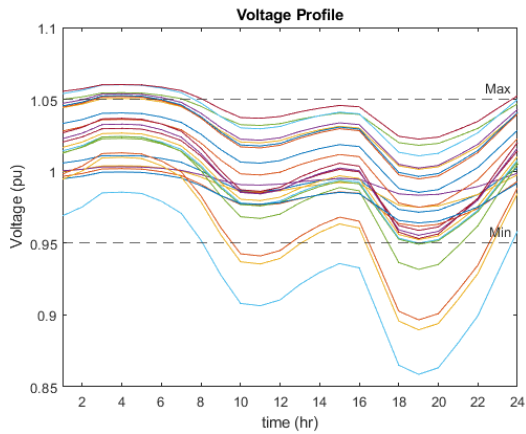

Figure 9. IEEE 30-bus system: Voltage profile (no controller)

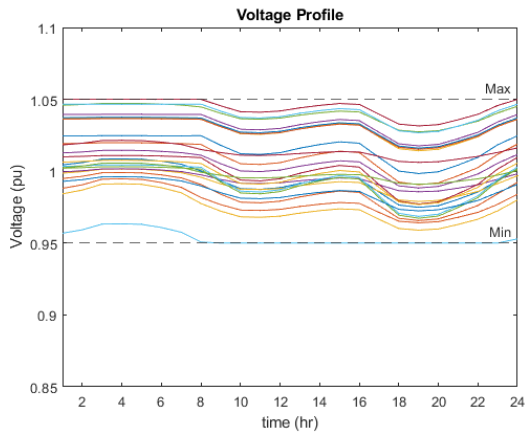

Figure 10. IEEE 30-bus system: Voltage profile (with controller)

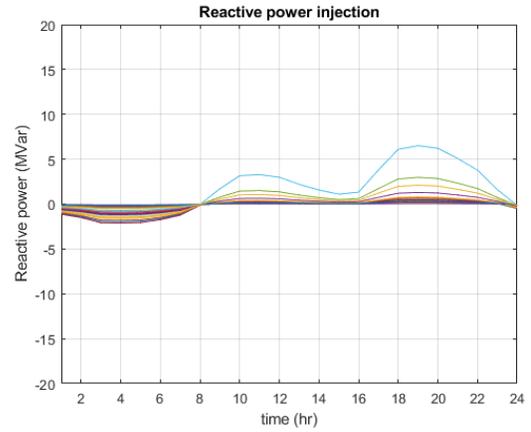

Figure 11. IEEE 30-bus system: Reactive power injection.

## IV. CONCLUSION

In this paper, optimal distributed feedback voltage controller was proposed. The performance was tested on two IEEE bus systems under static load and time varying load with time span of one day and 1 hour resolution. The controller managed to keep the voltage profile within acceptable magnitudes while satisfying the reactive power constrains of the control components in the optimum way in terms of operational cost. For future work, including real power to the control scheme is our goal.